\documentclass[12pt, a4paper]{article}
\usepackage{global}
\usepackage{booktabs}
\usepackage{tikz}
\usepackage{calrsfs}

\newcommand{\XX}{\mathbb{X}}
\newcommand{\LL}{\mathbb{L}}
\DeclareMathOperator{\Lip}{Lip}
\newcommand{\holder}[2]{\mathcal C^#1(#2)}

\newcommand{\h}{\mathfrak h}
\newcommand{\hh}{{h\vee\mathfrak\h}}
\newcommand{\Hn}{\mathcal H_n}

\newcommand{\g}{\bar g}
\newcommand{\gp}{\tilde{g}}

\newcommand{\x}{{x_0}}
\newcommand{\ra}{\mathfrak a}
\newcommand{\rb}{\mathfrak b}
\newcommand{\rc}{\mathfrak c}
\newcommand{\kd}{\mathfrak d}
\newcommand{\qq}{\mathfrak q}

\newcommand{\hf}{\hat f}

\newcommand{\finf}{\mathfrak B}
\newcommand{\re}{\mathfrak C}

\newcommand{\lipk}{L}

\newcommand{\lipholder}{\mathfrak L}

\definecolor{violet}{rgb}{0.8,0,0.8}

\let\E=\e
\let\P=\p

\usepackage{fourier-orns}

\begin{document}
\title{Pointwise Adaptive Estimation of the Marginal Density of
  a Weakly Dependent Process}
\author{Karine Bertin\thanks{CIMFAV, Universidad de
  Valpara\'{i}so, General Cruz 222,
  Valpara\'{i}so, Chile, tel/fax: 0056322303623}
\and
Nicolas Klutchnikoff\thanks{Crest-Ensai and  Université de Strasbourg}}
\date{\today}
\maketitle

\begin{abstract}
    This paper is devoted to the estimation of the common marginal density function of weakly dependent processes.
    The accuracy of estimation is measured using pointwise risks. We propose a data-driven procedure using kernel rules. The bandwidth is selected using the approach of Goldenshluger and Lepski and we prove that the resulting estimator satisfies an oracle type inequality. The procedure is also proved to be adaptive (in a minimax framework) over a scale of Hölder balls for several types of dependence: stong mixing processes, $\lambda$-dependent processes or \iid sequences can be considered using a single procedure of estimation. Some simulations illustrate the performance of the proposed method.
    {%
    }
\end{abstract}


\section{Introduction}\label{sec:intro}

Let $\XX=(X_i)_{i\in\Z}$ be a real-valued weakly dependent process admitting a common marginal density $f:\R\to\R$. We consider the problem of estimating $f$ at a fixed point $\x$ based on observation of $X_1,\ldots, X_n$ with $n\in\N^*$. The accuracy of an estimator $\tilde f_n$
is evaluated using the pointwise risk
defined, for fixed  $\x\in\R$ and $q>0$, by
\begin{equation}
    R_q(\tilde f_n, f) = \left(\e\abs{\tilde f_n(\x)-f(\x)}^q\right)^{1/q},
\end{equation}
where $\e$ denotes the expectation with respect to the distribution of the process $\XX$. The main interest in considering such risks is to obtain estimators that adapt to the local behavior of the density function to be estimated.

The aim of this paper is to obtain adaptive estimators of $f$ on H\"older classes of regularity $s>0$ for this risk and different types of weakly dependent processes.

In the independent and identically distributed (\iid) case, the minimax rate of convergence is $n^{-s/(2s+1)}$ (see \citet{tsybakov2004}, and references therein). Adaptive procedures based on the classical Lepki procedure (see \citet{lepski-1990}) have been obtained (see \citet{MR1772224}) with rates of the form $\left(\frac{\log n}{n}\right)^{s/(2s+1)}$.

In the context of dependent data, \citet{RW2006} as well as \citet{Rio} studied kernel density estimators from a minimax point of view for pointwise risks. They obtained the same rate of convergence as in the independent case when the coefficients of dependence decay at a geometric rate.
Several papers deal with the adaptive estimation of the common marginal density of a {weakly} dependent process. \citet{TV1998} and \citet{CM2002} proposed $\LL_p$--adaptive estimators under $\alpha$--mixing or $\beta$--mixing conditions that converge at the previously mentioned rates. 
\citet{GW2010} extend these results to a wide variety of {weakly} dependent processes including $\lambda$--dependent processes. 
Note that, in these papers, the proposed procedures are based on nonlinear wavelet estimators and only integrated risks are considered. Moreover, the thresholds are not explicitly defined since they depend on an unknown multiplicative constant. As a consequence, such methods can not be used directly for practical purposes.

Our main purpose is to prove similar results for pointwise risks.
We propose here a kernel density estimator with a data-driven selection of the bandwidth where the selection rule is performed using the so-called Goldenshluger-Lepski method \citep[see][]{MR2543590,GL2011,MR3230001}. This method was successfully used in different contexts such as in \citet{comteGenon2012,DHRR2012,BertinLacourRivoirard2014,rebelles2015}, but only with \iid observations. However there are at least two practical motivations to consider dependent data. Firstly, obtaining estimators that are robust with respect to  slight perturbations from the \iid\ ideal model can be useful. Secondly, many econometric models
(such as ARCH or GARCH)
deal with dependent data that admit a common marginal density. These two motivations suggest to consider a class of dependent data as large as possible and to find a single procedure of estimation that adapts to each situation of dependence.



Our contribution is the following. We obtain the adaptive rate of convergence for pointwise risks over a large scale of H\"older spaces in several situations of dependence, such as  $\alpha$--mixing introduced by \citet{R1956} and  the $\lambda$--dependence defined by \citet{DW2007}. This partially generalizes previous results obtained in i.i.d. case by \citet{MR1772224} and \citet{rebelles2015}. To the best of our knowledge, this is the first adaptive result for pointwise density estimation in the context of dependent data. To establish it, we prove an oracle type inequality: the selected estimator performs almost as well as the best estimator in a given large finite family of kernel estimators.
Our data-driven procedure depends only on explicit quantities. This implies that this procedure can be directly implemented in practice. As a direct consequence, we get a new method to choose an accurate local bandwidth for kernel estimators.

The rest of this paper is organized as follows. Section~\ref{sec:model} is devoted to the presentation of our model and of the assumptions on the process $\XX$. The construction of our procedure of estimation is developed in Section~\ref{sec:procedure}. The main results of the paper are stated in Section~\ref{sec:results} whereas their proofs are postponed to Section~\ref{sec:proofs}. A simulation study is performed in Section~\ref{sec:simulations} to illustrate the performance of our method in comparison to other classical estimation procedures. The proofs of the technical results are presented in the appendix.

\section{Model}\label{sec:model}

In what follows $\XX=(X_i)_{i\in\Z}$ is a real-valued discrete time process and the observation consists of the {vector $(X_1,\ldots,X_n)$.}
We assume that the $X_i$'s are identically distributed and we aim at estimating the common marginal density $f$ at a fixed point $x_0\in\R$. {In this section, basic assumptions on the distribution of $\XX$ are stated. Examples of processes are given, which illustrate the variety of models covered in this work.

\subsection{Assumptions}

We first assume that $f$ is bounded in a neigborhood of $\x$.
\begin{assumption}\label{ass:f}
    The marginal density $f$ satisfies
\begin{equation*}
 \sup_{x\in V_n(\x)}f(x)\leq \finf,
\end{equation*}
where $V_n(\x)=[\x-\left(\frac{1}{\log n}\right)^{1/q},\x+\left(\frac{1}{\log n}\right)^{1/q}]$ and $\finf$ is a positive real constant.

\end{assumption}

Before giving the two last assumptions, we need some notation. For $i\in\Z^u$, we consider the random variable $\XX_i=(X_{i_1},\ldots,X_{i_u})$ with values in $\R^u$.  If $i\in\Z^u$ and $j\in\Z^v$, we define the gap between $j$ and $i$  by $\gamma(j,i)=\min(j)-\max(i)\in\Z$.
For any positive integer $u$, the functional class $\mathbb{G}_{u}$  consists of real-valued functions $g$ defined on $\R^u$ such that the support of $g$ is included in $\left(V_n(\x)\right)^u$,
\begin{equation}
	\norme[\infty,u]{g}=\sup_{x\in[-1,1]^u}|g(x)|<\infty,
\end{equation}
and
\begin{equation*}
  \Lip_u(g) = \sup_{x\neq y} \frac{|g(x)-g(y)|}{\sum_{i=1}^u|x_i-y_i|}<+\infty.
\end{equation*}
We now define the sequence $\rho(\XX) =(\rho_r(\XX))_{r\in\N^*}$ by
\begin{equation*}
  \rho_r(\XX) =
  \sup_{u,v\in\N^*}\sup_{\substack{(i,j)\in\Z^u\times\Z^v \\ \gamma(j,i)\geq r}}
  \sup_{g\in\mathbb G_{u}}\sup_{\gp\in\mathbb G_{v}}\frac{\big|\cov\big(g(\XX_i),\gp(\XX_j)\big)\big|}{\Psi(u,v,g,\gp)},
\end{equation*}
where $\Psi(u,v,g,\gp)=\max(\Psi_1(u,v,g,\gp),\Psi_2(u,v,g,\gp))$ where for $g\in\mathbb{G}_u$ and $\gp\in\mathbb{G}_v$
\begin{equation*}
 	\Psi_1(u,v,g,\gp)= 4\norme[\infty,u]{g}\norme[\infty,v]{\gp}
 \end{equation*}  and
 \begin{equation*}
 	\Psi_2(u,v,g,\gp) = u\norme[\infty,v]{\gp}\Lip_{u}(g)+v\norme[\infty,u]{g}\Lip_{v}(\gp)+uv\Lip_u(g)\Lip_v(\gp).
 \end{equation*}

Both Assumptions \ref{ass:rho} and \ref{ass:cov} given below, allow us to control the covariance of functionals of the process. Assumption \ref{ass:rho} deals with the type of dependence, whereas Assumption \ref{ass:cov} is technical (see Lemma~\ref{lem:cov}, Section~\ref{sec:proofs}).

\begin{assumption}\label{ass:rho}
  For some positive constants $\ra$, $\rb$, and $\rc$, the sequence $\rho(\XX)$ is such that
  \begin{equation*}
  \rho_r(\XX) \leq \rc\exp(-\ra r^\rb), \quad \forall r\in\N.
 \end{equation*}
\end{assumption}

\begin{assumption}\label{ass:cov}
  The exists a positive constant $\re$ such that
  \begin{equation*}
  	\mu(\XX) = \sup_{g,\gp\in\mathbb{G}_1\setminus\{0\}}\sup_{i\neq j} \frac{\left|\e\left(g(X_i)\gp(X_j)\right)\right|}{\norme[1]{g}\norme[1]{\gp}} \leq \re,
  \end{equation*}
  where $\norme[1]{g}=\int_{\R}|g(t)|dt$.
\end{assumption}

\subsection{Comments}\label{sec:comments}

\paragraph{On Assumption~\ref{ass:f}.}
We are assuming that the marginal density $f$ is bounded on a neighborhood of $\x$. Such an assumption is classical in density estimation \citep[see][and references therein]{GL2011}. Note also that stationarity of $\mathbb X$ is not assumed. Thus, re-sampled processes of stationary processes can be considered in this study as in \citet{RW2006}.

\paragraph{On Assumption~\ref{ass:rho}.}
Recall that a process $\XX$ is called weakly dependent if, roughly speaking, the covariance between functionals of the past and the future of the process decreases as the gap from the past to the future increases (for a broader picture of weakly-dependent processes, as well as examples and applications, we refer the reader to \citet{MR2338725} and references therein). Assumption~\ref{ass:rho} ensures that the decay occurs at a geometric rate. Under similar assumptions, \citet{DN2007,MPR2009} proved Bernstein-type inequalities that are used in the proof of Theorem~\ref{thm:1}.

Note also that the type of weak-dependence considered in this paper includes usual types of dependence such as strong mixing as well as classical weak-dependence assumptions used in econometrics as illustrated below.

\bigskip
{\itshape \textbf{Strongly mixing processes.}}
The process $\mathbb{X}$ is called $\alpha$-mixing if the sequence $\alpha(\mathbb{X})=(\alpha_r(\mathbb{X}))_{r\in\N}$ defined by
\begin{equation*}
  \alpha_r(\mathbb{X}) = \sup_{n\in\Z} \sup_{A\in\mathcal F_{-\infty}^n}\sup_{B \in\mathcal F_{n+r}^{+\infty}}\abs{\p(A\cap B)-\p(A)\p(B)}
\end{equation*}
tends to $0$ as $r$ goes to infinity, where $\mathcal F_{k}^\ell$ is defined for any $k,\ell\in\Z$ as the $\sigma$-algebra generated by $(X_i)_{k\leq i\leq\ell}$. Recall that for any $u,v\in\N^*$, $g:\R^u\to\R$ and $\gp:\R^v\to\R$ such that $\norme[\infty,u]{g}<+\infty$ and $\norme[\infty,v]{\gp}<+\infty$ we have
\begin{equation}\label{eq:covmixing}
   \sup_{\substack{(i,j)\in\Z^u\times\Z^v \\ \gamma(j,i)\geq r}}\abs{\Cov(g(\mathbb{X}_i),\gp(\XX_j))}\leq 4 \norme[\infty,u]{g}\norme[\infty,v]{\gp} \alpha_r(\mathbb{X}).
\end{equation}
This readily implies that, for any $r\in\N$, we have $\rho_r(\XX)\leq \alpha_r(\XX)$.


\bigskip
{\itshape \textbf{$\lambda$--dependent processes.}}
The process $\XX$ is called $\lambda$--dependent if the sequence $(\lambda_r(\XX))_r$ defined by
\begin{equation*}
  \lambda_r(\XX) =
  \sup_{u,v\in\N^*}\sup_{\substack{(i,j)\in\Z^u\times\Z^v \\ \gamma(j,i)\geq r}}
  \sup_{g\in{\mathbb G}_{u}}\sup_{\gp\in{\mathbb G}_{v}}\frac{\big|\cov\big(g(\XX_i),\gp(\XX_t)\big)\big|}{\Psi_2(u,v,g,\gp)},
\end{equation*}
tends to $0$ as $r$ tends to infinity. Then  $\rho_r(\XX)\leq\lambda_r(\XX)$ for any $r\in\N$.

{\itshape Example.} Bernoulli shifts are, under some conditions, $\lambda$--dependent. Indeed, let us consider the process $\mathbb X$ defined by:
 \begin{equation}\label{eq:bernoulli-shift}
     X_i = H((\xi_{i-j})_{j\in\Z}), \quad i\in\Z
 \end{equation}
 where $H:\R^\Z\to[0,1]$ is a measurable function and the variables $\xi_i$ are \iid\ and real-valued. In addition, assume that there exists a sequence $(\theta_r)_{r\in\N^*}$ such that
 \begin{equation}\label{eq:theta}
     \e\abs{H((\xi_i)_{i\in\Z})-H((\xi_i')_{i\in\Z})}\leq \theta_r
 \end{equation}
 where, for any $r\in\N^*$, $(\xi_i')_{i\in\Z}$ is an \iid\ sequence such that $\xi'_i=\xi_i$ if $\abs{i}\leq r$ and $\xi'_i$ is an independent copy of  $\xi_i$ otherwise.
 It can be proved \citep[see][]{DL1999} that such processes are strongly stationary and $\lambda$--dependent with rate $\rho_r \leq 2\theta_{[r/2]}$.
 Remark also that $\theta_r$ can be evaluated under both regularity conditions on the function $H$ and integrability conditions on the $\xi_i$, $i\in\Z$. Indeed, if we assume that there exist $b\in(0,1]$ and positive constants $(a_i)_{i\in\Z}$ such that $\abs{H((x_i)_i)-H((y_i)_i)}\leq \sum_{i\in\Z}a_i\abs{x_i-y_i}^b$ with $\xi_i\in\LL_b(\R)$ for all $i\in\Z$, then
 \begin{equation}
     \theta_r = \sum_{\abs{i}\geq r} a_i \e\abs{\xi_i}^b.
 \end{equation}

 Moreover, under the weaker condition that $(\xi_i)_{i\in\Z}$ is $\lambda$--dependent and stronger assumptions on $H$ \citep[see][]{DW2007}, the process $\mathbb X$ inherits the same properties. Finally, we point out that classical econometrics models such as AR, ARCH or GARCH can be viewed as causal Bernoulli shifts (that is, they obey \eqref{eq:bernoulli-shift} with $j\in\N$).

%
%

\paragraph{On Assumption~\ref{ass:cov}}
This technical assumption is satisfied in several situations. In what follows, we offer examples of sufficient conditions such that Assumption~\ref{ass:cov} holds.

\bigskip
{\itshape Situation 1.} We assume that for any $i,j\in\Z$ the pair $(X_i,X_j)$ admits a density function $f_{i,j}$ with respect to the Lebesgue measure in $\R^2$. Moreover we assume that there exists a constant $\mathfrak F$ such that
\begin{equation*}
  \sup_{i,j\in\Z}\sup_{x,y\in V_n(\x)} f_{i,j}(x) \leq \mathfrak F.
\end{equation*}
Under this assumption and using Fubini's theorem, we readily obtain that  Assumption~\ref{ass:cov} holds with $\re=\mathfrak F$.
Note that this assumption is used in \citet{GW2010}.

\bigskip
{\itshape Situation 2.} We consider the infinite moving average process, with \iid  innovations $(\xi_j)_{j\in\Z}$, given by :
\begin{equation*}
  X_i = \sum_{j\in\Z} a_{j} \xi_{i-j},\qquad t\in\Z,
\end{equation*}
where $(a_j)_{j\in\N}$ and $(a_{-j})_{j\in\N}$  are decreasing sequences of deterministic positive real numbers. We assume that $\xi_1$ admits a density function $p_\xi(\cdot)$ bounded above by a positive constant $\mathfrak P$. Set $i,j\in\Z$ and denote by $\Xi$ the $\sigma$-algebra generated by the process $(\xi_i)_{i\in\Z\setminus\{i,j\}}$. For $g$, $\gp$ in $\mathbb{G}_1$ we have
\begin{align*}
    \e\left(g(X_i)\gp(X_j)\right)
    &= \e\left[\e\left(g(X_i)\gp(X_j)| \Xi\right)\right]\\
    &= \e\left[\e\left(g(a_0\xi_i+a_{j-i}\xi_j+A)\gp(a_0\xi_{j}+a_{i-j}\xi_i+B)| \Xi\right)\right],
\end{align*}
where $A$ and $B$ are $\Xi$-mesurable random variables. A simple change of variables gives:
\begin{equation*}
  \abs{\e\left(g(X_i)\gp(X_j)| \Xi\right)} \leq \frac{\mathfrak P^2}{a_0^2-a_{j-i}a_{i-j}}\norme[1]{g}\norme[1]{\gp}.
\end{equation*}
Since $a_0^2-a_{j-i}a_{i-j}\geq a_0^2-a_{1}a_{-1}$, Assumption~\ref{ass:cov} is fulfilled with $\mathfrak C = \mathfrak P^2(a_0^2-a_1a_{-1})^{-1}$.

\bigskip
{\itshape Situation 3.}
We consider a GARCH$(1,1)$ model. Let $\alpha$, $\beta$ and $\gamma$ be positive real numbers. Let $(\xi_i)_{i\in\Z}$ be \iid innovations with marginal density $p_\xi(\cdot)$, bounded above by a positive constant $\mathfrak{B}$, and denote by $(\mathcal F_i)_{i\in\Z}$ the natural filtration associated with this proccess. Assume that the process $\XX$ is such that, for any $i\in\Z$:
\begin{equation}\label{eq:GARCH}
  X_i = \sigma_i \xi_i
  \qquad\text{with}\qquad
  \sigma_i^2 = \gamma + \alpha X_{i-1}^2 + \beta \sigma_{i-1}^2.
\end{equation}
Consider $i,j\in\Z$ such that $i<j$. For $g$, $\gp$ in $\mathbb{G}_1$ we have
\begin{align}
    \e\left(g(X_i)\gp(X_j)\right)
    &= \e\e\left[g(X_i)\gp(X_j)|\mathcal F_i\right]\\
    &= \e\left(g(X_i)\e\left[\gp(X_j)|\mathcal F_i\right]\right).
\end{align}
Now remark that, since $\sigma_j\in\mathcal F_{j-1}$ and $\xi_j$ is independent of $\mathcal F_{j-1}$, we have
\begin{align}
    \e\left[\gp(X_j)|\mathcal F_i\right]
    &= \e\left(\e\left[\gp(\xi_{j}\sigma_j)|\mathcal F_{j-1}\right]| \mathcal F_i\right)\\
    &= \int_\R \gp(x) \e\left(\frac1{\sigma_{j}}p_\xi\left(\frac x{\sigma_{j}}\right)|\mathcal F_i\right)dx.
\end{align}
Since $\sigma_{j}\geq \sqrt{\gamma}$ we obtain:
\begin{equation}
  \abs{\e\left(g(X_i)\gp(X_j)\right)} \leq \frac{\mathfrak B}{\gamma} \norme[1]{g}\norme[1]{\gp}.
\end{equation}
Assumption~\ref{ass:cov} is thus fulfilled with $\mathfrak C =\mathfrak B/\gamma$.

\section{Estimation procedure}\label{sec:procedure}

In this section, we describe the construction of our procedure which is based on the so-called Goldenshluger-Lepski method (GLM for short). It consists in selecting, in a data driven way, an estimator in a given family of linear kernel density estimators. Consequently, our method offers a new approach to select an optimal bandwidth for kernel estimators in order to estimate the marginal density of a process in several situations of weak dependence. This leads to a procedure of estimation which is well adapted to inhomogeneous smoothness of the underlying marginal density. Notice also that our procedure is completely data-driven: it depends only on explicit constants that do not need to be calibrated by simulations or using the so-called rule of thumb.

\subsection{Kernel density estimators}\label{sec:kernel}

We consider kernels $K:\R\to\R$ that satisfy the following assumptions.

\begin{assumption}\label{ass:k1}
    The kernel $K$ is compactly supported on $[-1,1]$. Its Lipschitz constant $\Lip(K)$ is finite and $\int_\R K(x)\D x=1$.
\end{assumption}

\begin{assumption}\label{ass:k2}
    There exists $m\in\N$ such that the kernel $K$ is of order $m$. That is, for any $1\leq \ell\leq m$, we have
    {\begin{equation}
        \int_\R K(x) x^\ell \D x = 0 .   \end{equation}}
\end{assumption}
Let $h_* = n^{-1}\exp\big(\sqrt{\log n}\big)$ and $h^* = (\log n)^{-1/q}$ be two bandwidths and define $\Hn=\{2^{-k} : k\in\N\}\cap[h_*, h^*]$. We consider the family of estimators $\hat f_h$ defined by,
\begin{equation}
  \hat f_h(\x) = \frac1n\sum_{i=1}^n K_h(\x-X_i), \qquad h\in\mathcal H_n
\end{equation}
where $K_h(\cdot)=h^{-1} K(h^{-1}\cdot)$.

\subsection{Bandwidth selection}

Following \citet{GL2011}, we first define for $h\in\Hn$ the two following quantities
\begin{equation}
    A(h,\x) = \max_{\h\in\Hn}\left\{\abs{\hat f_{\hh}(\x)-\hat f_{\h}(\x)}-\widehat{M}_n(h,\h)\right\}_+
\end{equation}
and
\begin{equation}\label{eq:MNH}
 \widehat M_n(h) = \sqrt{2q\abs{\log h}\left(\widehat J_n(h)+\frac{\delta_n}{nh}\right)},
\end{equation}
where $\delta_n=(\log n)^{-1/2}$, $\{y\}_+=\max(0,y)$ for any $y\in\R$,
and $\hh=\max(h,\h)$ for any $h,\h\in\Hn$. We also consider
\begin{equation}\label{eq:MNHH}
 \widehat M_n(h,\h) = \widehat M_n(\h)+\widehat M_n(h\vee\h),
\end{equation}
 and
\begin{equation}
\widehat J_n(h)=\frac{1}{n^2}\sum_{i=1}^n K_h^2(x_0-X_i).
\end{equation}
Then our procedure consists in selecting  the bandwidth $\hat h(\x)$ such that
\begin{equation}\label{eq:selection-rule}
    \hat h(\x) = \argmin_{h\in\Hn} \left(A(h,\x)+\widehat M_n(h)\right).
\end{equation}
The final estimator of $f(\x)$ is defined by
\begin{equation}
    \hat f(\x) = \hat f_{\hat h(\x)}(\x).
\end{equation}

\begin{remark} The Goldenshluger-Lepski method consists in selecting a data-driven bandwidth that makes a trade-off between the two quantities $A(h,\x)$ and $\widehat M_n(h)$. Hereafter we explain how the minimization in \eqref{eq:selection-rule} can be viewed as an empirical version of the classical trade-off between a bias term and a penalized standard deviation term.
\begin{enumerate}
\item The quantity $\widehat{M}_n(h)$ can be viewed as a penalized upper bound of the standard deviation of the estimator $\hat{f}_h$. Indeed, Lemma~\ref{lem:stochastic-term} implies that
\begin{equation*}
   \Var(\hat f_h(\x))
   \leq J_n(h)+\frac{\delta_n}{6nh}
\end{equation*}
where
\begin{equation}\label{eq:Jn}
J_n(h)=\frac{1}{n}\int K_h^2(x_0-x)f(x)dx
\end{equation}
would be the variance of $\hat f_h(\x)$ if the data were \iid
Moreover (see the proof of Theorem~\ref{thm:1} in Section~\ref{sec:proofs}), for $n$ large enough and with high probability
\begin{equation*}
   J_n(h)+\frac{\delta_n}{6nh}\le \hat{J}_n(h)+\frac{\delta_n}{nh}.
\end{equation*}

\item The quantity  $A(h,\x)$ is a rough estimator of the bias term of $\hat{f}_h$. Indeed (see proof of Theorem~\ref{thm:1}), we have
\begin{align}
A(h,\x)\le& \max_{\h\in\mathcal H_n}|\e\hat f_{\hh}(\x) - \e\hat f_{\h}(\x)|+2 T \nonumber\\
\le & 2\max_{\h\leq h}\abs{K_\h\star f(\x)-f(\x)}+2T,\label{eq:bias}
\end{align}
where
\begin{equation*}
T=\max_{h\in\Hn} \left\{\abs{\hat f_{h}(\x)-\e\hat f_{h}(\x)}-\widehat{M}_n(h)\right\}_+.
\end{equation*}
The quantity $T$ is negligible with respect to $1/\sqrt{nh}$ with high probability
and $\max_{\h\leq h}\abs{K_\h\star f(\x)-f(\x)}$ is of the same order of the bias of $\hat{f}_h$ over H\"older balls.
\end{enumerate}
\end{remark}

\section{Results}\label{sec:results}


We prove two results. Theorem~\ref{thm:1} is an oracle-type inequality: under appropriate assumptions, our estimator performs almost as well as the best linear kernel estimator in the considered family. Theorem~\ref{thm:2} proves that our procedure achieves classical minimax rates of convergence (up to a multiplicative logarithmic factor) over a wide scale of Hölder spaces. 


\begin{theorem}\label{thm:1}
	Under Assumptions~\ref{ass:f},~\ref{ass:rho},~\ref{ass:cov} and~\ref{ass:k1} we have:
    \begin{align}\label{eq:thm1}
        R_q^q(\hat f,f)
        &\leq
        C^*_1   \min_{h\in\Hn} \left\{\max_{\substack{\h\leq h\\\h\in\Hn}}\abs{K_\h \star f(x_0)-f(x_0)}^q+\left(\frac{\abs{\log h}}{nh}\right)^{q/2}\right\}
    \end{align}
    where $C^*_1$ is a positive constant that depends only on $\ra$, $\rb$, $\rc$, $\finf$, $\re$ and $K$.
\end{theorem}

Proof of Theorem~\ref{thm:1} is postponed to Section~\ref{sec:proofs}.

\begin{remark}
     The right hand side term of \eqref{eq:thm1} can be viewed as a tight upper bound for $\min_{h\in\Hn} \e|\hat
    f_h(\x)-f(\x)|^q$ since it is the sum of an approximation of the bias term and the standard deviation term (up to a multiplicative logarithmic term) of $\hat f_h$. That means that our
    procedure performs almost as well as the best kernel density estimator in
    the considered family.
\end{remark}
Now using Theorem~\ref{thm:1}, we obtain in Theorem~\ref{thm:2} the adaptive rate of convergence on H\"older classes. Let $s$, $\mathfrak L$ and $\finf$ be positive real numbers. The H\"older class $\holder{s}{\lipholder,\finf}$ is defined as the set of functions $f:\R\to\R$ such that Assumption~\ref{ass:f} is fulfilled with the constant $\finf$, $f$ is $m_s=\sup\{k\in\N : k<s\}$ times differentiable and
\begin{equation}
    \abs{f^{(m_s)}(x)-f^{(m_s)}(y)} \leq \lipholder \abs{x-y}^{s-m_s},\qquad\forall x,y\in\R.
\end{equation}

\begin{theorem}\label{thm:2}
	Let $\ra$, $\rb$, $\rc$, $\finf$, $\re$, $m$  and $\mathfrak L$  be positive constants.
    Let $K$ be a kernel such that Assumptions~\ref{ass:k1} and \ref{ass:k2} are fulfilled (in particular, $K$ is a kernel of order $m$) and set $s$ such that $m_s\leq m$. There exists a constant $C_2^*$ that depends only on $\ra$, $\rb$, $\rc$, $\finf$, $\re$, $s$, $\mathfrak L$, $K$ and $q$ such that:
    \begin{equation}
        \sup_{\rho(\XX)\in \mathcal R(\ra,\rb,\rc)}\sup_{\mu(\XX)\leq\re}\sup_{f\in\holder{s}{\lipholder,\finf}}
        \left(\e\abs{\hat f(\x)-f(\x)}^q\right)^{1/q}
        \leq C_2^*\left(\frac{\log n}{n}\right)^{\frac{s}{2s+1}},
    \end{equation}
    where
    \begin{equation}
    	\mathcal R(\ra,\rb,\rc) = \left\{(\rho_r)_{r\in\N} : \rho_r \leq \rc \exp(-\ra r^{\rb})\right\}.
    \end{equation}
\end{theorem}
This result is a direct consequence of Theorem~\ref{thm:1}, since it can be easily proved that
\begin{align}
   \sup_{f\in\holder{s}{\lipholder,\finf}}
 \max_{\substack{\h\leq h\\\h\in\Hn}}\abs{K_\h \star f(x_0)-f(x_0)}^q
    \leq C_3^* h^{sq},
\end{align}
for any bandwidth $h>0$, where $C_3^*$ depends only on $s$, $\mathfrak L$, $K$ and $q$. This implies that, for $n$ large enough, there exists $h_n(s,\mathfrak L, K, q)\in\mathcal H_n$ such that the right hand side of \eqref{eq:thm1} is bounded, up to a multiplicative constant, by the expected rate.

\begin{remark}
	\begin{enumerate}
\item Recall that the expectation $\E$ is taken with respect to the distribution of the process $\XX$. Note also that the sequence $\rho(\XX)$, $\mu(\XX)$ and $f$ depend only on this distribution. As a consequence our procedure of estimation is minimax (up to a multiplicative $\log n$ term) with respect to any distribution of $\XX$ that satisfies the conditions:
\begin{equation}
	\rho(\XX)\in \mathcal R(\ra,\rb,\rc),
	\quad
	\mu(\XX)\leq\re
	\quad\text{and}\quad
	f\in\holder{s}{\lipholder,\finf}.
\end{equation}
Indeed, in the i.i.d. case (which is included in our framework since $\rho(\XX)\equiv 0$ and $\mu(\XX)\leq \finf^2$), the minimax rate of convergence over the H\"older class $\holder{s}{\mathfrak L,\finf}$ is of order $n^{-s/(2s+1)}$ and can be obtained from the results of \citet{MR1062695} or \citet{tsybakov2004}.
Moreover, note that $\hat f$ does not depend on the constants $\ra$, $\rb$, $\rc$, $\finf$, $\re$, $\mathfrak L$ and $s$ that appear in these conditions. Thus, our procedure is adaptive, up to the $\log n$ term, to both the regularity of $f$ and the ``structure'' of dependence.
\item
It can be deduced from our proofs that the minimax rate of convergence over H\"older classes, under Assumptions~\ref{ass:f}, \ref{ass:rho} and \ref{ass:cov}, is upper bounded, up to a multiplicative constant, by $n^{-s/(2s+1)}$. This result was previously obtained in a similar setting by  \citet{RW2006}  and \citet{Rio}. Given that this rate is minimax optimal in the \iid case,  it is also the minimax rate of convergence under our assumptions.

\item The extra $\log n$ term in the rate of convergence obtained in Theorem~\ref{thm:2} is unavoidable. Indeed, for pointwise estimation, even in the \iid case \citep[see][among others]{lepski-1990,klutchnikoff2014,rebelles2015} the adaptive rate of convergence is of this form. This ensures that our procedure attains the adaptive rate of convergence over H\"older classes.
 \item Note that $\delta_n$, that appears in \eqref{eq:MNH}, allows us to control the covariance terms of the development of $\Var(\hat f_h(\x))$ under Assumption~\ref{ass:rho}. If we only consider the \iid case, the covariance terms vanish, and the choice $\delta_n=0$ can be considered. The resulting procedure still satisfies an oracle inequality and remains adaptive in this case.
\item As far as we know, this result is the first theoretical pointwise adaptive result for the estimation of the marginal density in a context of weak dependence. Moreover, integrating the pointwise risk on a bounded domain, we obtain that our procedure converges adaptively at the rate $\left({n}^{-1}{\log n}\right)^{{s}/{(2s+1)}}$ in $\mathbb L_p$--norm ($p\neq \infty$) over H\"older balls. This extends the results of \citet{GW2010}.
\end{enumerate}
\end{remark}

\section{Simulation study}\label{sec:simulations}

In this section, we study the performance of our procedure using simulated data.  More precisely, we aim at estimating three density functions, for three types of dependent processes. In each situation, we study the accuracy of our procedure as well as classical competitive methods, calculating empirical risks using $p=10000$ Monte-Carlo replications.  In the following, we detail our simulation scheme and comment the obtained results.\\

\noindent
\textbf{Simulation scheme}\\

\noindent\textit{Density functions.} We consider three density functions to be estimated. The first one is:
\begin{equation}
f_1(x)=1.28\left(\sin((3\pi/2-1)x)I_{[0,0.65]}(x)+I_{(0.65,1]}(x)+cI_{[0,1]}(x)\right),
\end{equation}
where $c$ is a positive constant such that $f_1$ is a density.
The second one is the density of a mixture of three normal distributions restricted to the support $[0,1]$
\begin{equation}
f_2(x)=\left(\frac{1}{2}\phi_{0.5,0.1}(x)+\frac{1}{4}\phi_{0.6,0.01}(x)+\frac{1}{4}\phi_{0.65,0.95}(x)+c\right)I_{[0,1]}(x),
\end{equation}
where $\phi_{\mu,\sigma}$ stands for the density of a normal distribution with mean $\mu$ and standard deviation $\sigma$ and $c$ is a positive constant such that $f_2$ is a density.
Note that very similar densities were also considered in \citet{GW2010}.
The third one is:

\begin{equation*}
f_3(x) =\sum_{k=1}^{5} \left(2- 40\left|x-\frac{k}{10}+\frac{1}{20}\right|\right)I_{(\frac{k-1}{10},\frac{k}{10}]}(x)  +0.5I_{(0.5,1]}(x).
\end{equation*}

The function $f_1$ is very smooth except in the discontinuity point $x=0.65$.
The function $f_2$ is a classical example where rule-of-thumb bandwidths do not work. The third function has several spikes in $[0,0.5]$ and is constant on $[0.5,1]$. As a consequence, a global choice of bandwidth can fail to catch the two different behaviors of the function.
The three densities are bounded from  above (Assumption~\ref{ass:f} is then satisfied) and their inverse cumulative distribution functions are Lipschitz.\\

\noindent\textit{Types of dependence:} We simulate data $(X_1,\ldots,X_n)$ with density $f\in\{f_1,f_2,f_3\}$ in three cases of dependence. Denote by $F$ the cumulative distribution function of $f$.

\begin{enumerate}
\item[] \textbf{Case 1.} The $X_i$ are independent variables given by $F^{-1}(U_i)$ where the $U_i$ are i.i.d. uniform variables on $[0,1]$. Assumptions~\ref{ass:rho} and \ref{ass:cov} are clearly satisfied.
\item[] \textbf{Case 2.} The $X_i$ are $\lambda$--dependent given by $F^{-1}(G(Y_i))$ where the $Y_i$ satisfy of the non-causal equation:
\begin{equation}
  Y_i = 2(Y_{i-1}+Y_{i+1})/5+5\xi_i/21,\quad i\in\Z.
\end{equation}
Here $(\xi_i)_{i\in\Z}$ is an \iid sequence of Bernoulli variables with parameter $1/2$. The function $G$ is the marginal distribution function of the $Y_i$ and of the variable $\frac{U+U'+\xi_0}{3}$, where $U$ and $U'$ are independent uniform variables on $[0,1]$ \citep[see][for more details]{GW2010}. 

\item[] \textbf{Case 3.} The $X_i$ are given by $F^{-1}(G(Y_i))$  where $\mathbb{Y}=(Y_i)_{ i\in\Z}$ is an ARCH(1) process given by  \eqref{eq:GARCH} where the  $(\xi)_{i\in\Z}$ are \iid standard normal  variables, $\alpha=0.5$, $\beta=0$ and $\gamma=0.5$. In this case the function $G$ is estimated using the empirical distribution function on a simulated process $\mathbb{\tilde{Y}}$ independent of $\mathbb{Y}$ with the same distribution.
\end{enumerate}

\noindent
It remains to verify that Assumptions~\ref{ass:rho} and~\ref{ass:cov} hold for the processes $\mathbb{X}$ in the last two cases. Firstly, note that the process $\mathbb{Y}$ is $\lambda$--dependent with
\begin{equation}\label{lambda}\lambda_r(\mathbb{Y})\le \mathfrak{c}\exp(-\mathfrak{a}r),
\end{equation}
for some positive constants $\mathfrak{a}$ and $\mathfrak{c}$.  Indeed, in the second case, $Y_i$ is of the form \eqref{eq:bernoulli-shift} since it satisfies $Y_i=\sum_{j\in\Z} a_j\xi_{i-j}$ with $a_j=\frac{1}{3}\left(\frac{1}{2}\right)^{|j|}$  and the sequence $\theta_r$ (see \eqref{eq:theta}) such that $\theta_r\propto \left(\frac{1}{2}\right)^{r}$. In the third case, since $\alpha<1$, $Y_i$ is $\beta$--mixing at a geometric rate and then it is $\alpha$--mixing and $\lambda$--dependent at a geometric rate.

Now,  in both cases, since $F^{-1}$ and $G$ are Lipschitz and the process $(Y_i)_{i\in\Z}$ is $\lambda$--dependent, using Proposition 2.1 of \citet{MR2338725}, we have that  the process $(X_i)_{i\in\Z}$ is also $\lambda$--dependent  with $\lambda_r(\mathbb{X})$ satisfying \eqref{lambda} with some positive constants $\mathfrak{a}$ and $\mathfrak{c}$. As a consequence, Assumption~\ref{ass:rho} is fulfilled in the last two cases.

Secondly, Assumption~\ref{ass:cov} holds for the process $(Y_i)_{i\in\Z}$ (see \textit{Situation 2} and \textit{Situation 3} in Subsection~\ref{sec:comments}).
Then, for $g,\gp\in\mathbb{G}_1$ and $i,j\in \mathbb{Z}$ with $i\neq j$, we have
\begin{align*}
\left|\e\left(g(X_i)\gp(X_j)\right)\right|\le &\re\|g\circ F^{-1}\circ G\|_1 \|\gp\circ F^{-1}\circ G\|_1\\
\le &  \frac{\re \finf^2}{\mathfrak D^2} \|g\|_1\|\gp\|_1
\end{align*}
where $\mathfrak D=\min\{G'(x):x\in G^{-1}\circ F(V_n(\x))\}$ is bounded from below as soon as $\x\in(0,1)$  and $n$ is large enough. As a consequence, Assumption~\ref{ass:cov} is also fulfilled for the procces $\XX$.\\

\noindent
\textbf{Estimation procedures}\\

We propose to compare in this simulation study the following procedures.
\begin{itemize}
\item Our procedure (GL) $\hat{f}$ performed with, in \eqref{eq:MNH}, $q=2$ and $\delta_n=(\log n)^{-1/2}$.
\item The leave-one-out cross validation (CV) performed on the family of kernel estimators given in Subsection~\ref{sec:kernel}.
\item The kernel procedure with bandwidth given by the rule-of-thumb (RT).
\item The procedure developed by \citet{GW2010}.
\end{itemize}
In the first three procedures, we use the uniform kernel.\\

\noindent
\textbf{Quality criteria}\\

For each density function $f\in\{f_1,f_2,f_3\}$ and each case of dependence, we simulate $p=10000$ sequences of observations $(X_1,\ldots,X_n)$ with $n=1000$. Given an estimation procedure, we calculate $p$ estimators $\hat f^{(1)},\ldots,\hat f^{(p)}$.
We consider the empirical integrated square error:

\begin{equation*}
\widehat{ISE}=\frac{1}{p}\sum_{j=1}^p\int_{[0,1]}\left(f(x)-f^{(j)}(x)\right)^2 dx.
\end{equation*}

\noindent
\textbf{Results}\\

Our results are summarized in Table \ref{table:MISE}.

\begin{table}
\begin{center}
\begin{tabular}{|r|r|r|r|r|}
  \toprule
&  & Case 1 & Case 2 & Case 3\\
  \midrule
 & GL & 0.036 & 0.033 &  0.044 \\
 $f_1$ & CV & 0.027 & 0.034 & 0.049 \\
 & RT & 0.036 & 0.040 & 0.054 \\
 \hline
 & GL & 0.181 & 0.203 & 0.222 \\
 $f_2$ & CV & 0.079 & 0.116 & 0.162 \\
 & RT & 0.965 & 0.975 & 0.971 \\
 \hline
 & GL & 0.090 & 0.098 & 0.118 \\
 $f_3$ & CV & 0.172 & 0.180 & 0.190 \\
 & RT & 0.263 & 0.266 & 0.286 \\

  \bottomrule
\end{tabular}
\caption{Mean of ISE for the two densities $f_1$ and $f_2$, the three cases of dependence and the three procedures GL, CV and RT.}
\label{table:MISE}
\end{center}
\end{table}

For the estimation of the function $f_1$, our procedure gives better results than the CV or RT methods in cases of dependence (2 and 3).  We also outperform the results of \citet{GW2010} for case 1 and 2  where the ISE was around 0.09 (case 3 was not considered in their paper).
For the estimation of $f_2$ which is quite smooth, the cross validation method is about two times better than GL method and as expected, the RT method does not work.
For the estimation of $f_3$ that contains several peaks, the GL procedure is about two times better than the CV method.


To conclude, in the considered examples, our procedure has similar or better performances than already existing methods used for dependent data. Moreover, it gives better results when the density function to be estimated presents irregularities. This illustrates the fact that our method adapts locally to the irregularities of the function thanks to the use of local bandwidths. An other important point is that the choice of the bandwidth depends on explicit constants that can be used directly in practice and do not need previous calibration. Additionally, our GL procedure is about 25 times faster than cross-validation.

\section{Proofs}\label{sec:proofs}

\subsection{Basic notation}

For the sake of readability, we introduce in this section some conventions and {notations} that are used throughout the proofs. Moreover, here and later, we assume that Assumptions~\ref{ass:f}, \ref{ass:rho} and \ref{ass:cov} hold.

{Firstly,} let us consider, for any $h\in\Hn$, the functions $g_h$ and $\bar g_h$ defined, for any $y\in\R$, by $g_h(y) =  K_h(\x-y)$ and
\begin{equation}
  \bar g_h(y) = \frac{g_h(y)-\e g_h(X_1)}{n}.
\end{equation}
Note that we have:
\begin{align}
  \hf_h(\x) -\E\hf_h(\x) &= \sum_{i=1}^n \bar g_h(X_i),\quad h\in\Hn.
\end{align}
Next, we introduce some constants. Let us consider:
\begin{equation}
  C_1 = \norme[1]K^2(\finf^2+\re),
  \qquad
  C_2 = \finf\norme[2]K^2,
  \qquad\text{and}\qquad
  C_3 = 2\norme[\infty]K.
\end{equation}
Moreover we define $\lipk = \Lip(K)$ and
\begin{equation}
  C_4 = 2C_3\lipk+\lipk^2
  \quad
  \qquad\text{and}\qquad
  C_5=(2C_1)^{3/4}C_4^{1/4}.
\end{equation}
\subsection{Auxiliary results}

In the three following lemmas, assume that Assumptions~\ref{ass:f}, \ref{ass:rho} and \ref{ass:cov}  are satisfied.
The first lemma provides bounds on covariance terms for functionals of the past and the future of the observations. The considered functionals depend on the kernel $K$.
\begin{lemma}\label{lem:cov}
  For any $h\in\Hn$, we define
  \begin{equation}
    D_1(h) = D_1(n,h) = \frac{C_3}{nh}
    \qquad\text{and}\qquad
    D_2(h) = D_2(n,h) = \frac{C_5}{n^2h}.
  \end{equation}
  Then for any $u$, $v$ and $r$ in $\N$, if $(i_1,\ldots,i_u,j_1,\ldots j_v)\in\Z^{u+v}$ is such that $i_1\leq\ldots,i_u\leq i_u+r\leq j_1\leq\ldots\leq j_v$, we have
  \begin{equation}
    \abs*{\cov\left(\prod_{k=1}^u\bar g_h(X_{i_k}), \prod_{m=1}^v\bar g_h(X_{j_m})\right)}\leq \Phi(u,v)D_1^{u+v-2}(h)D_2(h) \rho_r^{1/4},
  \end{equation}
  where $\Phi(u,v)=u+v+uv$.
\end{lemma}

The following lemma provides a moment inequality for the classical kernel estimator.
\begin{lemma}\label{lem:stochastic-term}
There exists a positive integer $N_0$ that depends only on $\ra$, $\rb$, $\rc$, $\finf$, $\re$ and $K$ such that, for any $n\geq N_0$, we have
\begin{equation}
     \E\left(\left| \sum_{i=1}^n \g_h(X_i)\right|^{2}\right) \leq J_n(h)+\frac{\delta_n}{6nh}\leq \frac{C_2}{nh} + \frac{\delta_n}{6nh},
  \end{equation}
  where $J_n$ is defined by \eqref{eq:Jn}. Moreover for $\qq>0$, we have
  \begin{equation}
     \E\left(\left| \sum_{i=1}^n \g_h(X_i)\right|^{\qq}\right) \leq \mathfrak C_\qq(nh)^{-\qq/2}(1+o(1)),
  \end{equation}
  where $\mathfrak C_{\qq}$ is a positive constant. Here the $o(\cdot)$--terms depend only on $\ra$, $\rb$, $\rc$, $\finf$, $\re$ and $K$.
\end{lemma}

The following result is an adaptation of the Bernstein-type inequality obtained by \citet{DN2007}.
\begin{lemma}[Bernstein's inequality]
\label{lem:Bernstein}
We have:
\begin{equation}
  \P\left(\left|\sum_{i=1}^n \g_h(X_i)\right|\geq \lambda(t)\right) \leq C_7\exp(-t/2)
\end{equation}
where,
\begin{equation}\label{eq:lambda}
  \lambda(t) = \sigma_n(h)\sqrt{2t} + B_n(h)\left(2t\right)^{2+1/\rb},\qquad t\geq 0
\end{equation}
\begin{equation}\label{eq:sigman}
  \sigma_n(h) = J_n(h)+\frac{\delta_n}{6nh}
\end{equation}
and \begin{equation}\label{eq:Bn}
  B_n(h)= \frac{C_6}{nh}
\end{equation}
with $C_6$ and $C_7$ positive constants that depend only on $\ra$, $\rb$, $\rc$, $\finf$, $\re$ and $K$.
\end{lemma}

\subsection{Proof of Theorem~\ref{thm:1}}

Let us denote $\gamma=q(1+\delta_n/(12\max(C_2,1/6)))$. For convenience, we split the proof into several steps.

\noindent
\textbf{Step 1.} Let us consider the random event
\begin{equation}
\mathcal A = \bigcap_{h\in\Hn}\left\{  \left|\widehat J_n(h) - J_n(h)\right|\le \frac{\delta_n}{2nh}\right\}
\end{equation}
and the quantities $\Gamma_1$ and {$\Gamma_2$} defined by :
\begin{equation}
  \Gamma_1 = \E \left|\hat f(\x)-f(\x)\right|^q\1_{\mathcal A}
\end{equation}
and
\begin{equation}
  \Gamma_2 = \left(\max_{h\in\Hn}R_{2q}^{2q}(\hat f_h,f)\P(\mathcal A^c)\right)^{1/2}
\end{equation}
where $\1_{\mathcal A}$ is the characteristic function of the set $\mathcal A$. Using Cauchy-Schwarz inequality, it follows that:
\begin{align}
  R_q^q(\hat f, f) &\leq \Gamma_1 + \Gamma_2\label{eq:step1}.
\end{align}
We define
\begin{equation}
    \mathfrak M_n(h,a) =\sqrt{2q\abs{\log h}\left( J_n(h)+\frac{a\delta_n}{nh}\right)}.
\end{equation}
Now note that if the event $\mathcal A$ holds, we have:
\begin{equation}\label{eq:mnh}
    \mathfrak M_n\left(h,\frac{1}{2}\right) \leq \widehat M_n(h) \leq \mathfrak M_n\left(h,\frac{3}{2}\right).
\end{equation}
Steps 2--5 are devoted to control the term $\Gamma_1$ whereas $\Gamma_2$ is upper bounded in Step 6.

\noindent
\textbf{Step 2.} Let $h\in\mathcal{H}_n$ be an arbitrary bandwidth. Using triangular inequality we have:
\begin{equation}
|\hat{f}(x_0)-f(x_0)|\le |\hat{f}_{\hat{h}}(x_0)-\hat{f}_{h\vee\hat{h}}(x_0)|+|\hat{f}_{h\vee\hat{h}}(x_0)-\hat{f}_{h}(x_0)|+|\hat{f}_{h}(x_0)-f(x_0)|.
\end{equation}
If $h\vee\hat{h}=\hat{h}$, the first term of the right hand side of this inequality is equal to $0$, and if $h\vee\hat{h}=h$, it satisfies
\begin{align}
  |\hat{f}_{\hat{h}}(x_0)-\hat{f}_{h\vee\hat{h}}(x_0)|
  &\leq \left\{|\hat{f}_{\hat{h}}(x_0)-\hat{f}_{h\vee\hat{h}}(x_0)|-\widehat M_n(h,\hat h)\right\}_+ + \widehat M_n(h,\hat h)\\
  &\leq \max_{\h\in\Hn}\left\{|\hat{f}_{\h}(x_0)-\hat{f}_{\hh}(x_0)|-\widehat M_n(h,\h)\right\}_+ + \widehat M_n(h,\hat h)\\
  &\leq A(h,\x) + \widehat M_n(\hat h)+\widehat M_n(h).
\end{align}
Applying the same reasoning to the term $|\hat{f}_{h\vee\hat{h}}(x_0)-\hat{f}_{h}(x_0)|$ and
using \eqref{eq:selection-rule}, this leads to
\begin{equation}
|\hat{f}(x_0)-f(x_0)|\le 2\big(A(h,x_0)+\widehat M_n(h)\big)+|\hat{f}_{h}(x_0)-f(x_0)|.
\end{equation}
Using this equation, we obtain that, for some positive constant $c_q$,
\begin{equation}\label{eq:thm-1}
  \Gamma_1 \leq c_q \left(\E \big(A^q(h,\x)\1_{\mathcal A}\big) + \mathfrak M_n^q\left(h,\frac{3}{2}\right) + R_q^q(\hat f_h, f))\right).
\end{equation}

\noindent
\textbf{Step 3.}
Now, we upper bound $A(h,\x)$. Using basic inequalities we have:
\begin{align}
  A(h,\x) &\leq
  \max_{\h\in\Hn}
  \left\{
   \left|\E \hat{f}_{\hh}(x_0)-\E\hat{f}_{\h}(x_0)\right|
  \right\}_+ +
  \max_{\h\in\Hn}
  \left\{
   \left|\hat{f}_{\h}(x_0)-\E\hat{f}_{\h}(x_0)\right|
  -\widehat{M}_n(\h)\right\}_+\\
  &
  \qquad +
  \max_{\h\in\Hn}
  \left\{
   \left|\hat{f}_{\hh}(x_0)-\E \hat{f}_{\hh}(x_0)\right|
  -\widehat{M}_n(\hh)\right\}_+  \\
  &\leq \max_{\h\in\Hn}
  \left\{
   \left|\E \hat{f}_{\hh}(x_0)-\E\hat{f}_{\h}(x_0)\right|
  \right\}_++2 T,
\end{align}
where
\begin{equation}
        T=\max_{\h\in\Hn}
  \left\{
   \left|\hat{f}_{\h}(x_0)-\E\hat{f}_{\h}(x_0)\right|
  -\widehat{M}_n(\h)\right\}_+=\max_{\h\in\Hn}\left\{\left|\sum_{i=1}^n\g_\h(X_i)\right|-
            \widehat{M}_n(\h)\right\}_+.
\end{equation}
Using \eqref{eq:bias}, we obtain
\begin{align}
  A(h,\x) &\leq  2\max_{\substack{\h\leq h\\ \h\in\Hn }}\left\{\left|K_{\h}\star f(\x)-f(\x)\right|\right\}_+ + 2T\label{eq:thm-6}.
\end{align}
Denoting by $E_h(\x)$ the first term of the right hand side of \eqref{eq:thm-6}, we obtain
\begin{equation}\label{eq:step3}
  \E \left(A^q(h,\x)\1_{\mathcal A}\right) \leq c_q \left(E_h^q(\x)+\E \big(T^q\1_{\mathcal A}\big)\big)\right),
\end{equation}
for some positive constant $c_q$.

\noindent
\textbf{Step 4.}
It remains to upper bound $\E (T^q\1_{\mathcal A})$. To this aim, notice that,
\begin{equation}\label{eq:t-ttilde}
  \E (T^q\1_{\mathcal A}) \leq \E \tilde T^q,
\end{equation}
where
\begin{equation}
  \tilde T = \max_{\h\in\Hn}\left\{\left|\sum_{i=1}^n\g_\h(X_i)\right|-\mathfrak M_n\left(\h,\frac{1}{2}\right)\right\}_+.
\end{equation}}
Now, we define $r(\cdot)$ by
\begin{equation}
  r(u)=\sqrt{2\sigma_n^2(h)u} +
  2^{\kd-1}B_n(h)\left(2u\right)^{\kd}, \qquad u\geq 0
\end{equation}
where $B_n(h)$ and $\sigma_n(h)$ are given by~\eqref{eq:sigman}  and \eqref{eq:Bn} and $\kd=2+\rb^{-1}$ .
Since $h\geq h_*=n^{-1}\exp(\sqrt{\log n})$, we have, for $n$ large enough:
\begin{equation}\label{eq:toto1}
  2^{\kd-1}B_n(h)(2\gamma \abs{\log h})^{\kd} \leq \frac{\delta_n}{12\sqrt{2\max(C_2,1/3)}} \sqrt{\frac{2q\abs{\log h}}{nh}}.
\end{equation}
Moreover, we have
\begin{align}
\sqrt{2\sigma_n^2(h)\gamma\abs{\log h}} &\le \sqrt{2q\abs{\log h}\left(J_n(h)+\frac{\delta_n}{6nh}\right)\left(1+\frac{\delta_n}{12\max(C_2,1/6)}\right)}\nonumber\\
    &\le \sqrt{2q\abs{\log h}\left(J_n(h)+\frac{\delta_n}{3nh}\right)}.\label{eq:toto2}
\end{align}
Last inequality comes from the fact that $J_n(h)$ is upper bounded by ${C_2}/{(nh)}$. Now, using \eqref{eq:toto1} and \eqref{eq:toto2}, we obtain:
\begin{align}
  r(\gamma|\log h|)  &\leq \mathfrak M_n\left(h,\frac{1}{2}\right). \label{eq:penr}
\end{align}

Thus, doing the change of variables
$t=\left(r(u)\right)^{q}$ and thanks to~\eqref{eq:penr}, we obtain:
\begin{align*}
\E \tilde T^q &\le \sum_{\h\in\Hn}\int_0^{\infty} \P
\left(\left|\sum_{i=1}^n\g_\h(X_i)\right|\ge \mathfrak M_n\left(\h,\frac{1}{2}\right)+t^{1/q}\right)dt\\
&\le C\sum_{\h\in\Hn}\int_0^{\infty}
r'(u) r(u)^{q-1}\P
\left(\left|\sum_{i=1}^n\g_\h(X_i)\right|\ge r(\gamma|\log
  \h|)+r(u)\right)du,\\
&\le  C\sum_{\h\in\Hn}\int_0^{\infty} u^{-1} \lambda(u)^{q}\P \left(\left|\sum_{i=1}^n\g_\h(X_i)\right|\ge \lambda(\gamma|\log \h|+u) \right)du,
\end{align*}
where $\lambda(\cdot)$ is defined by~\eqref{eq:lambda}. Using Lemma~\ref{lem:Bernstein}, we obtain
\begin{align*}
  \E \tilde T^q &\le
  C\sum_{\h\in\Hn}\int_0^{\infty}u^{-1}\left(\sqrt{\sigma_n^2(\h)u}+B_n(\h)u^{3}\right)^{q}\exp\left\{-\frac{u}{2}-
    \frac{\gamma|\log \h|}{2}\right\}du\\
&\le
C\sum_{\h\in\Hn}\sigma_n^{q}(\h)\h^{\gamma/2}.
\end{align*}
Using the definitions of the quantities that appear in this equation and using \eqref{eq:t-ttilde}, we readily obtain:
\begin{align}
  \E T^q\1_{\mathcal A} \leq Cn^{-q/2}\sum_{k\in\N} (2^{(\gamma-q)/2})^{-k}
  &\leq C \frac{n^{-q/2}}{1-\exp\left(-\frac{\delta_n\log2}{24\max(C_2,1/6)}\right)}\\
  & \leq C n^{-q/2}\sqrt{\log n}\\
  & \leq C (nh)^{-q/2}.\label{eq:thm-2}
\end{align}

\noindent
\textbf{Step 5.}
Lemma~\ref{lem:stochastic-term} implies that:
\begin{equation}\label{eq:thm-4}
  \E|\hat f_h(\x)-f(\x)|^q \leq c_q\left(|E_h(\x)|^q+(nh)^{-q/2}\right)
\end{equation}
for some positive constant $c_q$.

Combining \eqref{eq:thm-1}, \eqref{eq:step3}, \eqref{eq:thm-2} and \eqref{eq:thm-4}, we have:
\begin{equation}\label{eq:control-gamma1}
  \Gamma_1 \leq C^* \min_{h\in\Hn} \left\{\max_{\substack{\h\leq h\\\h\in\Hn}}\abs{K_\h \star f(x_0)-f(x_0)}^q+\left(\frac{\abs{\log h}}{nh}\right)^{q/2}\right\}
\end{equation}
where $C^*$ is a positive constant that depends only on $\ra$, $\rb$, $\rc$, $\finf$, $\re$ and $K$.

\noindent
\textbf{Step 6.} Using Lemma~\ref{lem:Bernstein} where in $\overline{g}_h$, $K$ is replaced by $K^2$, we obtain that
\begin{equation*}
\P\left(\left|\frac{h}{n}\sum_{i=1}^nK_h^2(x_0-X_i)-h\int K^2_h(x_0-x)f(x)dx\right|> \frac{\delta_nx}{2} \right)\le \exp(-C_1n^2h^2),
\end{equation*}
where $C_1$ is a constant that depends only on $\ra$, $\rb$, $\rc$, $\finf$, $\re$ and $\delta$.
Then this implies that
\begin{equation*}
\P\left(\mathcal A^c\right)\le \frac{\log n}{\log 2}\exp\big(-C_1\exp(2\sqrt{\log n})\big)
\end{equation*}
and then
\begin{equation}\label{eq:thm-5}
\Gamma_2=o\left(\frac{1}{n^{q/2}}\right).
\end{equation}
Now, using~\eqref{eq:control-gamma1}  and
\eqref{eq:thm-5}, Theorem~\ref{thm:1} follows.

\appendix

\section{Proof of Lemma~\ref{lem:cov}}

In order to prove this lemma, we derive two different bounds for the
term
\begin{align}
  \Upsilon_h(u,v) =  \left|\cov \left( \prod_{k=1}^u \g_h(X_{i_k}), \prod_{m=1}^v
      \g_h(X_{j_m})\right)\right|.
\end{align}
The first bound is obtained by a direct calculation whereas the second
one is obtained thanks to the dependence structure of the
observations.
For the sake of readability, we denote $\ell=u+v$
throughout this proof.

\noindent
\textbf{Direct bound.}
The proof of this bound is composed of two steps. First, we assume that $\ell=2$, then the general case $\ell\geq 3$ is considered.

Assume that $\ell=2$.
If Assumptions~\ref{ass:f} and \ref{ass:cov} are fulfilled, we have
\begin{align}
  n^2\Upsilon_h(u,v) &\le  \left\vert\e\left(g_h(X_i)
      g_h(X_j)\right)\right\vert +\left(\e g_h(X_1)\right)^2\\
&\leq  {(\re+\finf^2)}\norme[1]{g_h}^2\leq C_1.
\end{align}
Then, we have
\begin{equation}\label{eq:cov-l-2}
  \abs*{\Cov\left(\g_h(X_i),\g_h(X_j)\right)} \leq C_1n^{-2}.
\end{equation}

Let us now assume that $\ell\geq 3$.
Without loss of generality, we
can assume that $u\geq2$ and $v\geq 1$. We have:
\begin{align}
  \Upsilon_h(u,v) &\leq A+B
\end{align}
where
\begin{align}
  \begin{cases}
    A = \E\left( \prod_{k=1}^u \g_h(X_{i_k}) \prod_{m=1}^v
      \g_h(X_{j_m})\right)\\
    B = \E\left( \prod_{k=1}^u \g_h(X_{i_k})\right)
    \E\left(\prod_{m=1}^v \g_h(X_{j_m})\right).\\
  \end{cases}
\end{align}
Remark that both $A$ and $B$ can be bounded, using~\eqref{eq:cov-l-2}, by
\begin{align}
  \norme[\infty]{\bar g_h}^{(u-2)+v}\Cov(\bar g_h(X_{i_1}),
  \bar g_h(X_{i_2})) &\leq \left(\frac{C_3}{nh}\right)^{\ell-2}\frac{C_1}{n^2}.
\end{align}
This implies our first bound, for all $\ell\geq 2$:
\begin{align}\label{eq:lem-2-first-part}
  \Upsilon_h(u,v) &\leq \frac{2C_1}{n^2}\left(\frac{C_3}{nh}\right)^{\ell-2}.
\end{align}

\noindent
\textbf{Structural bound.} Using Assumption~\ref{ass:rho}, we obtain that
\begin{align}
  \Upsilon_h(u,v) &\leq    \Psi\left(u,v,\g_h^{\otimes u},\g_h^{\otimes v}\right)
   \rho_r.
\end{align}
Now using that
\begin{equation*} \left\|  \frac{nh\bar g_h}{C_3}\right\|_\infty\le 1
\end{equation*}
and
\begin{equation*}
\Lip \left(\frac{nh\bar g_h}{C_3}\right)^{\otimes u} \le \Lip \left(\frac{nh\bar g_h}{C_3}\right)\le \frac{L}{C_3h},
\end{equation*}
we obtain, since $h\leq h^*$, that
\begin{align}
  \Upsilon_h(u,v) &\leq  \left(\frac{C_3}{nh}\right)^\ell \Phi(u,v)\frac{C_4}{C_3^2h^2}
   \rho_r.
\end{align}
This implies that
\begin{equation}\label{eq:struct2}
  \Upsilon_h(u,v) \leq  \frac1{n^2}\left(\frac{C_3}{nh}\right)^{\ell-2}\frac{C_4}{h^{4}}\Phi(u,v)\rho_{r}.
\end{equation}

\noindent
\textbf{Conclusion.}
Now combining \eqref{eq:lem-2-first-part} and~\eqref{eq:struct2} we
obtain:
\begin{align}
  \Upsilon_h(u,v) &\leq
  \frac{1}{n^2}\left(\frac{C_3}{nh}\right)^{\ell-2}
  (2C_1)^{3/4} \left(\frac{C_4}{h^4} \Phi(u,v) \rho_r
  \right)^{1/4}\\
  &\leq \frac{C_5}{n^2h}
  \left(\frac{C_3}{nh}\right)^{\ell-2} \Phi(u,v)\rho_r^{1/4}.
\end{align}
This proves Lemma~\ref{lem:cov}.
\section{Proof of Lemma~\ref{lem:stochastic-term}}
Proof of this result can be readily adapted from the proof of Theorem~1 in \citet{DL2001} (using similar arguments that ones used in the proof of Lemma~\ref{lem:cov}). The only thing to do is to bound explicitely the term
\begin{equation}
     A_2(\g_h) = \e \left(\sum_{i=1}^n \g_h(X_i)\right)^{\!2}.
 \end{equation}
Set $R=h^{-1/4}$. Remark that
\begin{align}
  A_2(\g_h) &= n\E\g_h(X_1)^2 +
  \sum_{i\neq j} \E\g_h(X_i)\g_h(X_j) \\
  &= J_n(h)+2\sum_{i=1}^{n-1}\sum_{r=1}^{n-i}\E\g_h(X_i)\g_h(X_{i+r}).
\end{align}
Using Lemma~\ref{lem:cov} and \eqref{eq:cov-l-2}, we obtain:
\begin{align}
  A_2(\g_h) &\leq  J_n(h)+2n\sum_{r=1}^R\frac{C_1}{n^2} +
  2D_2(h)\sum_{r=R+1}^{n-1}(n-r)\Phi(1,1)\rho_r^{1/4}\\
  &\leq
  J_n(h)+\frac{1}{nh}\left((2C_1)h^{3/4}+(6C_5)\sum_{r=R+1}^\infty\rho_r^{1/4}\right)\\
  &\leq
  J_n(h)+\frac{1}{nh}\left(\frac{2C_1}{(\log n)^{3/4}}+(6C_5)\sum_{r=1+(\log n)^{1/4}}^\infty\rho_r^{1/4}\right).
\end{align}
Last inequality holds sine $h\leq h^*\leq (\log n)^{-1}$. Using Assumption~\ref{ass:rho}, there exists $N_0 = N_0(K, \finf,\re, \ra,\rb,\rc)$ such that, for any $n\geq N_0$ we have:
\begin{equation}
    A_2(\g_h) \leq J_n(h)+\frac{\delta_n}{6nh}.
\end{equation}
This equation, combined with the fact that $J_n(h)\leq C_2(nh)^{-1}$, completes the proof.
\section{Proof of Lemma~\ref{lem:Bernstein}}

First, let us remark that Lemma~6.2 in \citet{GW2010} and Assumption~\ref{ass:rho} imply that there exist positive constants $L_1$ and $L_2$ (that depend on $\ra$, $\rb$ and $\rc$) such that, for any $k\in\N$ we have,
\begin{equation}
  \sum_{r\in\N}(1+r)^k\rho_r^{1/4} \leq L_1L_2^k(k!)^{1/\rb}.
\end{equation}
This implies that, using Lemma~\ref{lem:stochastic-term}, one can apply the Bernstein-type inequality obtained by \citet[see Theorem 1]{DN2007}. First, remark that, using Lemma~\ref{lem:stochastic-term}, for $n$ large enough, we have
\begin{equation}
  \e\left(\sum_{i=1}^n\bar g_h(X_i)\right)^{\!2} \leq \sigma_n(h)
  \quad\text{and}\quad
  B_n(h) = \frac{2L_2C_3}{nh}.
\end{equation}
where the theoretical expression of $B_n(h)$ given in \citet{DN2007}.
Let us now denote $\kd=2+\rb^{-1}$. We obtain:
\begin{align*}
  \P\left(\left|\sum_{i=1}^n \g_h(X_i)\right|\geq u\right) \leq
  \exp \left(-\frac{u^2/2}{\sigma_n^2(h) +
      B_n^{\frac1\kd}(h)  u^{\frac{2\kd-1}{\kd}}}\right).
\end{align*}
Now, let us remark that, on the one hand $\lambda(t)\geq \sigma_n(h)\sqrt{2t}$ and
thus $\lambda^2(t)\geq2\sigma_n^2(h)t$. On the other hand,
$\lambda^{2+\frac{1-2\kd}{\kd}}(t)\geq 2B_n^{\frac1\kd}(h)t$ and thus
$\lambda^{2}(t)\geq (2B_n^{\frac1\kd}(h)t)\lambda^{\frac{2\kd-1}{a\kd}}(t)$. This implies that
$\lambda^2(t)\geq t(\sigma_n^2(h)+B_n^{\frac1\kd}(h)\lambda^{\frac{2\kd-1}{\kd}}(t))$ and thus,
finally:
\begin{equation}
  \exp\left(-\frac{\lambda^2(t)}{\sigma_n^2(h)+B_n^{\frac1\kd}(h)\lambda^{\frac{2\kd-1}{\kd}}(t)}\right)\leq
  \exp(-t/2).
\end{equation}
This implies the results.

\section{Acknowledgments}
The authors have been supported by Fondecyt project 1141258. Karine
Bertin has been supported by the grant Anillo ACT-1112 CONICYT-PIA, Mathamsud project 16-MATH-03 SIDRE and ECOS project C15E05.

\bibliographystyle{plainnat}
\renewcommand*{\bibfont}{\small}
\bibliography{biblio}


\end{document}